%----------------------------------------------------------------
%------------------------- This is AMS-TeX   ---------------------
%-----------------------------------------------------------------
\input amstex
\documentstyle{amsppt}

\TagsOnRight

\nopagenumbers

\magnification=\magstep1

%Greek letters
                
\def\bet{\beta}         
                \def\Gam{\Gamma}
                
\def\eps{\varepsilon}           

\def\tet{\theta}                
\def\iot{\iota}

\def\lam{\lambda}

\def\ome{\omega}                \def\Ome{\Omega}

%Caligraphic roman letters

\def\calE{{\Cal E}}
\def\calF{{\Cal F}}

\def\calH{{\Cal H}}

\def\calL{{\Cal L}}
\def\calM{{\Cal M}}
\def\calN{{\Cal N}}

\def\calP{{\Cal P}}
\def\calQ{{\Cal Q}}

%Euler Fraktur letters

     \def\grg{{\frak g}}
     \def\grt{{\frak t}}

%Capital roman double letters 

\def\RR{\Bbb R}
\def\ZZ{\Bbb Z}
\def\CC{\Bbb C}

\def\otet{\overline\tet} 
\def\oxi{\overline\xi} 

\def\nek{,\ldots,}
\def\sdp{\times \hskip -0.3em {\raise 0.3ex
\hbox{$\scriptscriptstyle |$}}} % semidirect product

%words in roman font

\def\End{\operatorname{End\,}}
\def\const{\operatorname{const}}

\def\im{\operatorname {im}}
\def\ind{\operatorname{ind}}

%------------------------------------------------------------------------
%------------------------------------------------------------------------

\def\F{\calF}
\def\H{\calH}
\def\fo{\calF_{|_Z}\otimes o(Z)}

\def\E{\calE}

\def\w{\tet}

\def\rank{\operatorname{rank}}

\def\n{\nabla}
\def\nb{\nabla_T}

\def\Ss{\Cal S}

\def\au{Au}
\def\ab1{AB1}
\def\ab{AB2}
\def\bl{BL}
\def\bott{Bo1}
\def\bottb{Bo2}
\def\bfa{BF1}
\def\bfb{BF2}
\def\giv{G}
\def\far{Fa}
\def\hsi{Hs}
\def\mcd{McD}
\def\na{N1}
\def\nb{N2}
\def\nc{N3}
\def\pa{Pa}
\def\fr{Fr}
\def\O{O}
%-------------------------------------------------------------------
\topmatter

\title Equivariant Novikov inequalities
\endtitle
%\rightheadtext{}
%\leftheadtext{}
\author  Maxim Braverman and Michael Farber \endauthor
\address
School of Mathematical Sciences,
Tel-Aviv University,
Ramat-Aviv 69978, Israel
\endaddress
\email farber\@math.tau.ac.il, maxim\@math.tau.ac.il 
\endemail
\thanks{The research was  supported by grant No. 449/94-1 from the
Israel Academy of Sciences and Humanities}
\endthanks
\abstract{We establish an equivariant generalization of the Novikov 
inequalities which allow to estimate the topology of the set of critical 
points of a 
closed basic invariant form by means of twisted 
equivariant cohomology of the 
manifold. We apply these inequalities to study cohomology
of the fixed points set of a symplectic torus action. We show that in this case our
inequalities are {\it perfect}, i.e. they are in fact equalities.}
\endabstract
\endtopmatter
%-------------------------------------------------------------------------
%----------------------------------------------------------------------

In \cite{\na,\nb} S.P.Novikov associated to any real
cohomology class $\xi\in H^1(M,\RR)$ of a closed manifold $M$ 
a sequence of integers  \/
$\bet_0(\xi)\nek$ $\bet_n(\xi)$ (where $n=\dim M$) and then it was shown  
that for any closed
1-form $\w$ on $M$, having non-degenerate critical points, 
the Morse numbers $m_p(\w)$ satisfy the following inequalities
$$
	\sum_{i=0}^{p}(-1)^im_{p-i}(\w)\ 
	  \ge\  \sum_{i=0}^{p}
	    (-1)^i\beta_{p-i}(\xi), \qquad p=0,1,2,\dots,
$$
where $\xi=[\w]\ \in H^1(M,\RR)$ is the cohomology
class of $\w$.  Cf. also \cite{\far}, page 48, where stronger inequalities are 
mentioned.

In this paper we obtain an equivariant generalization of the Novikov
inequalities.  We consider a compact $G$ manifold $M$, where $G$ is a compact 
Lie group, and an invariant closed 1-form $\theta$ on $M$. Assuming that the
form $\theta$ is {\it basic} (cf. below) we define an equivariant generalization
of the Novikov numbers. To do so we show that any closed basic invariant 
form
determines a one parameter family of flat bundles over the Borel construction 
$M_G$, and the equivariant Novikov numbers are then defined as the 
{\it background value of the dimension of the cohomology} of 
this family. It turns out that these equivariant
Novikov numbers depend only on the cohomology class of the form $\theta$.
Using these equivariant Novikov numbers we construct {\it the Novikov
counting power series}, which is the first main ingredient appearing in 
our inequalities. 

We assume that the form $\theta$ {\it is non-degenerate in the sense of 
Bott}; this means that the critical points of $\theta$ form a submanifold $C$ 
such that the Hessian of $\theta$ is non-degenerate on the normal bundle 
to $C$.  Each
connected component $Z$ of the critical manifold $C$ determines the following
objects: an integer $\ind(Z)$ 
({\it the index}), a flat line bundle $o(Z)$ ({\it the orientation bundle})
and a subgroup $G_Z\subset G$ ({\it the stabilizer of the component}). Using 
these data we construct {\it the equivariant Morse counting power series} which
combines information about the equivariant cohomology of all the components
$Z$ of $C$. Our main theorem (Theorem 1.7) states, roughly, that 
the equivariant Morse power series {\it is greater} (in an appropriate sense) 
than the Novikov counting series. This statement gives actually an 
infinite number of 
inequalities involving the dimensions of the equivariant cohomology of
connected components of the critical set $C$ and the equivariant Novikov
numbers.

We use in this paper {\it equivariant cohomology twisted by an equivariant 
flat vector
bundle}, which is crucial for our approach. On one hand, any closed
invariant basic 1-form determines a one-parameter family of such bundles, 
which we 
use to define the equivariant generalizations of the Novikov numbers.
On the other hand we observe, that using this kind of cohomology allows
to strengthen the inequalities considerably. Very simple examples 
show (we were very surprised to discover this!) that applying the 
well-known equivariant Morse inequalities of Atiyah and Bott \cite{AB1, \bottb} 
to the 
case when 
the group $G$ is finite, 
one gets estimates, which
are sometimes worse than the standard Morse inequalities (ignoring 
the group action). We are going to devote a separate paper
describing in more detail what our approach gives in this classical situation; 
we will construct there a variety of new Morse type inequalities (labeled by 
the irreducible representations of the group $G$ ) and also some 
combinatorial analysis and comparison between them. These results are 
based on the results of the present paper.

As an application of the established in the present paper equivariant 
Novikov inequalities we obtain new relations between the Betti numbers of the 
set of fixed points of a symplectic torus action and the equivariant Novikov numbers.
It seems to us that 
the main novelty of this our result is that, in contrast with the
previously known bounds for the fixed points for symplectic torus actions 
using the equivariant cohomology
(cf. \cite{\giv}), our inequalities allow to find or to estimate the individual 
Betti numbers of
the fixed points set, and not only the sum of these numbers. We also show that 
under certain additional assumptions (for example, if the fixed points are all
isolated, cf. \S 5) it is possible to find an explicit 
expression for the numbers of fixed points of different indices 
in terms of the equivariant Novikov numbers. In the case of a holomorphic action
of the circle on a K\"ahler manifold preserving the K\"ahler form (which was
studied by E.Witten in \cite{Wi}) we obtain even simpler expressions involving
only the dimensions of the equivariant cohomology.

The proof of the main result of the paper (Theorem 1.7) is based in its 
main part
on the Novikov type inequalities for differential forms with
non-isolated zeros, obtained in \cite{\bfa,\bfb}.

The paper is organized as follows.  
In Section 1, we define the equivariant Novikov numbers, introduce necessary 
notations and then formulate our main result. We also discuss here some 
easy examples and special cases.
The sections 2 - 4 are devoted to the proof of the main theorem.
In Section 2, we recall some basic facts about equivariant vector
bundles, and twisted equivariant cohomology, which we use.
In section 3 we collected some remarks about basic 1-forms.
In Section 4, we prove the equivariant Novikov inequalities.
Finally, in Section 5 we describe our main application: the results about the 
Betti numbers of the
fixed points set of a torus action on a symplectic manifold.

Some results of this paper were briefly announced in \cite{BF3}. In \cite{BF3} we
developed also applications to the finite group case and to studying of
manifolds with boundary.  

We are thankful to J.Bernstein and L.Polterovich for useful discussions.
%---------------------------------------------------------------------
%---------------------------------------------------------------------

\heading \S 1. The main result \endheading

In this section we formulate the main result of the paper (Theorem 1.7).

Let $G$ be a compact Lie group and let $M$ be a closed
$G$-manifold. We do not demand that $G$ is connected and do not
exclude the case where $G$ is a finite group. We shall denote by
$\grg$ the Lie algebra of $G$.

%---------------------------------
\subheading{1.1. Basic 1-forms} The equivariant generalization of the Novikov 
inequalities which we will describe in this paper, are applicable to
closed 1-forms which are {\it basic}. Recall that a smooth 1-form $\w$ on 
a $G$ manifold $M$ is called basic (\cite{\ab}) if it is 
$G$-invariant and its restriction on any orbit of the action of $G$ equals to
zero. In other words, this means that the following two conditions have to be
satisfied
$$
	g^*\w=\w, \quad \iot(X_M)\w=0 \quad
		\text{for any} \quad g\in G, \ X\in \grg.\tag1
$$
Here $X_M$ is the vector field on $M$ defined by the infinitesimal
action of $\grg$ on $M$ and $\iot(X_M)$ denotes the interior
multiplication by $X_M$. 

Note that, if the group $G$ is finite, then $\w$ is basic if and only if 
it is $G$-invariant, i.e. if $g^*\w=\w$ is satisfied for any $g\in G$. 
Also, by Lemma 3.4, {\it if $M$ is connected and if the set of
fixed points of the action of $G$ on $M$ is not empty, then any closed 
$G$-invariant 1-form on $M$ is basic}. Note also that {\it any exact
invariant form $\theta=df$ is basic}.

Let $H_G^\ast(M,\RR)$ denote the equivariant cohomology of $M$. There
exists a natural map 
$$H_G^1(M,\RR)\to H^1(M,\RR)\tag2$$ 
(cf. \cite{\ab}). We observe as Lemma 3.3 that {\it a real 
cohomology class 
$\xi\in H^\ast(M,\RR)$ lies in the image of this map if and only if it 
may be represented by a basic differential form.}

%---------------------------------------------------
\subheading{1.2. Equivariant generalization of the Novikov numbers} In 
section 2 we discuss in detail the notion of {\it equivariant flat vector 
bundle over a $G$ manifold $M$}. We also construct {\it equivariant 
cohomology
$H^\ast_G(M,\F)$ twisted by an equivariant flat bundle} $\F$, cf. 2.7 - 2.9.
It is shown in 2.2 that {\it any closed basic 1-form $\theta$ determines
an equivariant flat line bundle $\E_\theta$} (with $\theta$ being its 
connection form). Using these constructions we define the equivariant
Novikov numbers as follows.

Given an equivariant flat bundle $\F$ over $M$ and a closed basic 1-form 
$\theta$ on $M$, consider the one-parameter family $\F\otimes \E_{t\theta}$
of equivariant flat bundles, where $t\in \RR$, ({\it the Novikov deformation}) 
and consider the twisted equivariant cohomology
$$H^i_G(M,\F\otimes\E_{t\theta}), \qquad\text{where}\quad t\in\RR,\tag3$$
as a function of $t\in \RR$.
The following Lemma describes the behavior of the dimension of the cohomology 
of this one-parameter family.

\proclaim{1.3. Lemma} For fixed $\F$ and $i$, there exists a finite
subset $S\subset\RR$ such that the dimension of the cohomology 
$H^i_G(M,\F\otimes\E_{t\theta})$ is constant for $t\notin S$ 
and the dimension of the cohomology $H^i_G(M,\F\otimes\E_{t\theta})$
jumps up for $t\in S$.
\endproclaim

The proof follows from Lemma 4.1 below if one takes into account 
Definition 2.9.

The subset $S$ which appears in the Lemma, is called the 
{\it set of jump points};
the value of the dimension of $H^i_G(M,\F\otimes\E_{t\theta})$ for 
$t\notin S$ is called {\it the background value of the dimension of this
family}.

\proclaim{1.4. Definition} The $i$-dimensional {\it equivariant Novikov 
number} $\bet_i^G(\xi,\F)$ is defined as the background value of the 
dimension of the cohomology of the family 
$H^i_G(M,\F\otimes\E_{t\theta})$.\endproclaim

Here $\xi\in \im[H_G^1(M,\RR)\to H^1(M,\RR)]$ denotes the cohomology class of    
$\theta$. By Lemma 2.3, the equivariant flat bundle $\E_\theta$ is determined
(up to gauge equivalence) only by $\xi$. By Lemma 3.3, any class 
$\xi\in \im[H_G^1(M,\RR)\to H^1(M,\RR)]$ can be realized by a basic closed    
form and so the equivariant Novikov numbers $\bet_i^G(\xi,\F)$ are
defined for all classes in the image $\xi\in\im[H_G^1(M,\RR)\to H^1(M,\RR)]$.

The formal power series
$$
	\calN_{\xi,\F}^G(\lam)=
		\sum_{i=0}^\infty\lam^i\bet_i^G(\xi,\F)\tag4 
$$ 
will be  called {\it the equivariant Novikov series}.

%-------------------------------------------------------------
\subheading{1.5} We will assume that the given basic 1-form 
$\w$ is {\it non-degenerate} in the sense of
R.Bott \cite{\bott}. Recall, that this means that the points of $M$, 
where the form $\w$ vanishes form a submanifold of $M$ (called the 
{\it critical set} $C$ of $\w$) and the {\it Hessian} of $\w$ is
non-degenerate on the normal bundle to $C$.

In order to make clear this definition, note that if we fix a tubular
neighborhood $N$ of $C$ in $M$, then the monodromy of $\w$ along
any loop in $N$ is obviously zero. Thus there exists a unique real
valued smooth function $f$ on $N$ such that $df=\w_{|_N}$ and
$f_{|_C}=0$. The {\it Hessian} of $\w$ is then defined as the
Hessian of $f$.

Let $\nu(C)$ denote the normal bundle of $C$ in $M$. Note that
$\nu(C)$ may have different dimension over different connected
components of $C$.  Since Hessian of $\w$ is non-degenerate, the
bundle $\nu(C)$ splits into the Whitney sum of two subbundles
$$
	\nu(C)\ =\ \nu^+(C)\ \oplus \nu^-(C),\tag5
$$
such that the Hessian is strictly positive on $\nu^+(C)$ and strictly
negative on $\nu^-(C)$. Here again, the dimension of the bundles
$\nu^+(C)$ and $\nu^-(C)$ over different connected components of the
critical point set may be different.

For every connected component $Z$ of the critical set $C$, the
dimension of the bundle $\nu^-(C)$ over $Z$ is called the {\it index}
of $Z$ (as a critical submanifold of $\w$) and is denoted by $\ind(Z)$.

Let $o(Z)$ denote the {\it orientation bundle} of $\nu^-(C)_{|_Z}$ 
considered as a flat real  line bundle.

%Let $G_0$ denote the connected component of $G$ which contains the unit
%element. Then $G_0$ is a normal subgroup of $G$. We denote by 
%$\Gam=G/G_0$ the quotient group. It is a finite group called the {em
%group of connected components} of $G$. 

%-----------------------------
\subheading{1.6. The equivariant Morse series}
Let $Z$ be a connected component of the critical set $C$.  If the
group $G$ is connected, then $Z$ is a $G$-invariant submanifold of
$M$. In the general case we denote by 
$$
	G_Z=\{g\in G| \ g\cdot Z\subset Z\}\tag6 
$$ 
the stabilizer of the component $Z$ in $G$.  Let $|G:G_Z|$ denote the
index of $G_Z$ as a subgroup of $G$.  Since $G_Z$ contains the
connected component of the unity in $G$, this index is finite.

The compact Lie group $G_Z$ acts on the manifold $Z$ and the flat
vector bundles $\F_{|_Z}$ and $o(Z)$ are $G_Z$-equivariant.  Let
$$
	H_{G_Z}^\ast(Z,\fo)\tag7
$$ 
denote the {\it equivariant cohomology} of the
flat $G_Z$-equivariant vector bundle $\fo$.  Consider the {\it
equivariant Poincar\'e series} of $Z$ 
$$
	\calP^{G_Z}_{Z,\F}(\lam)\ =\ 
	   \sum_{i=0}^{\infty} \lam^i \dim_{\CC}H_{G_Z}^i(Z,\fo)\tag8  
$$
and define using it the following {\it equivariant Morse counting
series}
$$
	\calM_{\w,\F}^G(\lam)\ = \ \sum_Z 
	     \lam^{\ind(Z)} |G:G_Z|^{-1}\calP_{Z,\F}^{G_Z}(\lam),\tag9
$$
where the sum is taken over all connected components $Z$ of $C$.

The main result of the paper is the following: 

\proclaim{1.7. Theorem} 
   Suppose that $G$ is a compact  Lie group, $\F$ is a flat
   $G$-equivariant vector bundle over a closed $G$-manifold $M$ and
   $\w$ is a closed non-degenerate (in the sense of Bott) basic 1-form
   on $M$. Then there exists a formal power series $\calQ(\lam)$ with
   non-negative integer coefficients, such that
   $$
	\calM_{\w,\F}^G(\lam)- \calN_{\xi,\F}^G(\lam)=
		(1+\lam)\calQ(\lam).   \tag10
   $$
\endproclaim
The theorem is proven in Section 4.

It is an interesting question under which conditions one may expect vanishing of the 
power series $\calQ(\lam)$ in (10). For the case $\xi=[\tet]=0$, this problem was solved by 
Atiyah and Bott \cite{AB1}. We show below in Theorem 5.3
that the generalized moment map in our symplectic
applications (cf. \S 5) {\it is perfect}, i.e. $\calQ(\lam)$ vanishes
in our inequalities (10).

We observe also, that our arguments work and prove
perfectness in a more general situation:
namely, under the conditions of Theorem 1.7, assuming
additionally that $G$ is connected and the set of critical points
of a basic 1-form (as in Theorem 1.7) coincides
with the fixed point set, cf. remark 5.4.

\subheading{1.8. Examples} Consider first the case when $G$ acts freely on 
$M$. Then the basic form      
$\theta$ defines a closed 1-form $\theta'$ on $M/G$ (cf. \cite{KN}, Ch.XII,
\S1). It is straightforward to 
see, that in this case the inequalities of Theorem 1.7 (with $\F$ being 
the trivial line bundle) reduce to the usual Novikov inequalities with 
respect to the form $\theta'$ on the quotient manifold $M/G$. In particular,
we see that in this situation the equivariant Novikov numbers 
$\beta_i^G(\xi, \F)$ vanish for large $i$.

As another example, consider the case of {\it the circle action} $G=S^1$. 
Applying the Localization Theorem (cf. 2.11) we obtain that 
for large $i$ the 
equivariant Novikov numbers $\beta_i^G(\xi, \F)$ are two-periodic and
coincide with the sum of even or odd, respectively,  usual 
(i.e. non-equivariant) Novikov numbers of the
fixed point set, cf. 2.12 and (56).

We shall consider now a special case of Theorem 1.7.

%---------------------------------------------------------------------
\subheading{1.9. The case of isolated critical points} In this section
we will assume that all critical points of $\w$ are isolated. For
simplicity we will assume that for any critical point $x\in C$ the
action of the stabilizer $G_x=\{g\in G:\ g\cdot x=x\}$ of $x$ on the
tensor product $\F_{|_x}\otimes o(\{x\})$ is trivial. Note that this
condition is automatically satisfied if the group $G$ is connected.

It is clear that any two singular points of $\theta$, belonging to the 
same orbit of $G$, have the same index. Let $m_i(\theta)$ denote {\it the 
number of orbits of critical points of $\theta$ having Morse
index $i$}. We will 
show that the equivariant Novikov inequalities established in Theorem 1.7 
give estimates of the numbers $m_i(\theta)$.

Denote by $H_{G}^\ast$ the $G$ equivariant cohomology ring of the point with 
complex coefficients, i.e. the cohomology of the classifying space
$BG$ of $G$ with coefficients in $\CC$. 
Set $\bet_i^G=\dim_{\CC} H^i_G$ and let
$$
	\calP^G(\lam)=\sum_{i=0}^\infty \lam^i\bet_i^G\tag11 
$$ 
be the Poincar\'e series of $G$. Note that, if $G$ is a finite group,
then $\calP^G(\lam)=1$. Also, if $H$ is a subgroup of $G$ of finite
index then $\calP^H(\lam)=\calP^G(\lam)$.

It is easy to see that the equivariant Morse counting polynomial (9) 
takes the form
$$
	\calM^G_{\w,\F}(\lam) = 
		d\cdot \calP^G(\lam)\sum_{i=0}^\infty \lam^i m_i(\w),
		\qquad\text{where}\quad d=\dim \F.  \tag12
$$
We observe, that the Morse counting polynomial depends in this case only 
on the dimension $d$ of the flat equivariant bundle $\F$. However, simple 
examples (constructed using remark 1.8) show that the Novikov counting 
polynomial $\calN_{\xi,\F}^G(\lam)$ 
may really depend on choice of the flat bundle $\F$. Theorem 1.7 gives 
in this situation many inequalities involving the numbers $m_i(\theta)$ 
and the equivariant Novikov numbers.

In the case when all critical points of $\theta$ are isolated, it is more realistic 
to assume that $G$ is a {\it finite group}. Then $\calP^G(\lam)=1$ and
we obtain
$\calM^G_{\w,\F}(\lam) = d\cdot \sum_{i=0}^\infty \lam^i m_i(\w)$. Thus 
Theorem 1.7 gives the inequalities
$$
	\sum_{i=0}^{p}(-1)^im_{p-i}(\w)\ \ge\  
		d^{-1}\cdot \sum_{i=0}^{p}(-1)^i\bet^G_{p-i}(\xi, \F), 
				\qquad p=0,1,2,\dots,   \tag13
$$
which look similar to the usual Novikov inequalities \cite{N1}, \cite{N2},
\cite{N3} (but note that the meaning of our notation $m_i(\theta)$ is
different).

\heading{\S 2. Twisted equivariant cohomology} \endheading

In this section we review some basic facts about equivariant flat vector
bundles, which we need in this paper. We are mainly interested in the  
pushforward construction which produces a flat vector bundle over the 
base space of a principal fiber bundle assuming that an equivariant flat
vector bundle over the space of the fibration is given. We will apply the
construction of pushforward in order to produce flat vector bundles
over the Borel construction $M_G$.

%-----------------------------------
\subheading{2.1. Equivariant flat vector bundles}
Let $G$ be a compact Lie group and let $M$ be a $G$-manifold. 
We assume that $G$ acts smoothly on $M$ from the left; we denote this action 
by

$$
	 G\times M \to M, \qquad (g,x)\mapsto g\cdot x.\tag14
$$

Let $\F\to M$ be a flat vector bundle over $M$ and let

$$\nabla:\Omega^\ast(M,\F)\to \Omega^{\ast+1}(M,\F),\qquad \nabla^2\ 
=\ 0\tag15$$ 
denote the covariant derivative defined by the flat structure on $\F$.
We will suppose that $G$ acts smoothly on the space $\F$ as well
$$
	G\times \F\to \F, \qquad (g,\xi)\mapsto g\cdot \xi 
$$ 
and that this action is compatible with the action of $G$ on $M$ in the 
sense that the projection of the bundle $\F\to M$ is $G$-equivariant.  
We will also assume
that the action $g: \F_x\to \F_{g\cdot x}$ is linear for any $g\in G$
and $x\in M$.

Now, the group $G$ acts on the space of smooth sections $\Omega^0(M,\F)$ 
by the formula 
$$
	(g\cdot s)(x)= g\cdot s(g^{-1}\cdot x).\tag16
$$
We will consider also the similar action of $G$ on the spaces of smooth 
forms $\Omega^\ast(M,\F)$ on $M$ with values in $\F$.

For any element $X$ of the Lie algebra $\grg$ of $G$, we will denote by
$\calL^\F(X)$ the corresponding infinitesimal action 
$$
	\calL^\F(X)s =\frac d{d\eps}{\big|_{\eps=0}}\exp(\eps X)\cdot s,
	\qquad\text{where}\quad s\in \Omega^0(M,\F).  \tag17
$$
Also, for $X\in \grg$, we denote by $X_M$ the vector field on $M$ defined by
the infinitesimal action of $\grg$ on $M$.

Now we arrive at the main definition: \newline
{\it A flat bundle $\F\to M$ as above, 
supplied with an action of $G$ on $\F$ which is linear on fibers and which
is compatible with the projection, is called  
$G$-equivariant flat vector
bundle if the following two conditions are satisfied:
$$g\circ \n\ =\ \n\circ g\  :\Omega^0(M,\F)\to \Omega^1(M,\F) \tag 18$$
for any $g\in G$, and
$$\n_{X_M}\ =\ \calL^\F(X):\ \Omega^0(M,\F)\to \Omega^0(M,\F) \tag 19$$
for any $X\in \grg$.} 

Note that condition (19) determines completely the covariant derivative 
in the 
directions tangent to the orbits of $G$ in terms of the action of $G$ on $M$
and on $\F$. It implies that {\it any flat section $s$ of $\F$ (determined 
locally over an open set $U\subset M$) is invariant, i.e. $s(gx)=gs(x)$ 
for all $x\in U$ and for any $g\in G$ close to the unit $e\in G$}.

In particular, it follows from the above remarks that condition (19) implies 
that the {\it flat bundle $\F$ is trivial over any orbit of $G$, if we 
assume that $G$ acts freely}.

Note also, that if we have two different flat connections $\nabla$ and 
$\nabla^\prime$ satisfying (19) then the difference $\nabla-\nabla^\prime$
is a 1-form with values in $\End(\F)$ which vanishes on vectors tangent to
the orbits of $G$, i.e. it is basic.

It is clear that if the group $G$ is
connected,  then condition (18) follows from (19). Conversely, 
if $G$ is a finite group, the second condition (19) carries no information.
Observe, that for a general compact group $G$ (which is neither connected,
nor finite) the conditions (18) and (19) are independent.    

%------------------------------------------------------------
\subheading{2.2. Example} The following example is crucial for our
study of the Novikov inequalities.
Let $M$ be a $G$ manifold and let $\theta$ be a 
closed $G$-invariant real 1-form on $M$. The last condition means that
$g^\ast\theta=\theta$ for all $g\in G$. Consider the flat vector bundle 
determined by the form $\theta$. Namely, let 
$$\E_\theta\ =\  M\times \CC\tag20$$ 
with the $G$ action coming from the factor $M$ and with the flat connection
$$\nabla \ =\ d\ +\ \theta\wedge\cdot\ \ .\tag21$$
Then this flat bundle always satisfy (18) and it satisfies (19) if and only 
if the form $\theta$ is basic.

\proclaim{2.3. Lemma} The equivariant flat bundle $\E_\theta$ (cf. above) 
considered up to equivariant gauge equivalence, depends only on the 
cohomology class $\xi=[\theta]\in H^1(M,\RR)$.
\endproclaim
\demo{Proof} We may assume that $M$ is connected; otherwise, one can 
consider the situation over the connected components separately.

Suppose that $\theta_1\ = \ \theta\ +\ df$ is another basic
1-form, where $f:M\to \RR$ is a smooth function. Then for
any $X\in\grg$ the derivative $X_M(f)$ vanishes and so the function $f$
must be constant on the connected components of the orbits of $G$. Since
$df$ is $G$-invariant we obtain that for any $g\in G$ there is a constant
$c(g)$ such that $f(gx)=f(x)+c(g)(x)$ for all $x\in M$. 

We claim that 
$c(g)$ must be identically zero. If there is $g\in G$ with $c(g)>0$ then 
we obtain $f(gx)>f(x)$ for all $x\in M$ which is impossible. To see this,
consider a point $x\in M$ such that $f$ achieves at $x$ it maximum on the 
compact set $Gx$. The case when $c(g)<0$ may be treated similarly.

The arguments above prove that the function $f$ is $G$ invariant. Now we 
may define the gauge equivalence $F:\E_{\theta_1}\to \E_\theta$ as the
operator of multiplication by the function $F=e^f$. It is equivariant,
commutes with the projection and satisfies 
$\nabla\circ F\ =\ F\circ\nabla_1$,
where $\nabla_1 \ =\ d\ +\ \theta_1\wedge\cdot\ $. $\square$
\enddemo

\subheading{2.4. Remarks} We first note, that given two equivariant flat 
vector bundles $\F_1$
and $\F_2$ over a $G$-manifold $M$, the Whitney sum $\F_1\oplus\F_2$ and
the tensor product $\F_1\otimes\F_2$ (considered together with their standard 
flat connections and the $G$ actions) are equivariant flat vector bundles.

Another easy observation: if $f:M_1\to M_2$ is an equivariant map between
$G$ manifolds and if $\F\to M_2$ is an equivariant flat bundle over $M_2$,
then the induced bundle $f^\ast\F$ over $M_1$ is also an equivariant 
flat bundle.

%------------------------------------------------------------
\subheading{2.5. The pushforward construction} We will
suppose here that {\it the action of $G$ on $M$ is free}. Then (by the Gleason
lemma, cf. \cite{\hsi}, page 9) the quotient space
$B=M/G$ is a smooth manifold and the quotient map $q:M\to B$ is a
locally trivial fibration. We will see that in this case any equivariant 
flat vector bundle over $M$ determines canonically a flat vector bundle
over $B$.

Let $\F$ be a $G$-equivariant flat vector bundle over $M$ and let
$\Ss(\F)$ denote the locally constant sheaf of flat sections of $\F$.
Denote by $q_*\Ss(\F)$ the direct image of sheaf $\Ss(\F)$. Since 
$\F$ is trivial over any connected component of any orbit of $G$ (cf. above), 
and since the projection $q$ is locally trivial,
it follows that $q_*\Ss(\F)$ is also a locally constant sheaf over $B$; 
its fiber over a point $b\in B$ can be identified with $\oplus_i \F_{m_i}$,
where $m_i$'s are some representatives of the connected components of
$q^{-1}(b)$.  Let $q_*\F$ denote the flat bundle corresponding to the
sheaf $q_*\Ss(\F)$.

Observe that the group $G$ acts naturally on the direct image $q_*\F$
and this action is compatible with the flat structure. The action of
the connected component $G_0$ of the unit element $e\in G$ is
obviously trivial. Thus, the flat bundle $q_*\F$ splits into a direct
sum of its flat subbundles corresponding to different irreducible
representations of the finite group $G/G_0$. The most important for us
will be the subbundle corresponding to the trivial representation; we
will denote it by
$$
	q^G_*(\F).\tag22
$$ 
It is the subbundle of $G$-invariants in $q_*\F$. 

In particular, if $G$ is
connected, then $q_*^G\F=q_*\F$.

{\it We will say that the flat bundle $q_*^G\F$ is the pushforward of
the bundle $\F$.} 
\subheading{\it Remark} We refer the reader to \cite{\bl} for a detailed
study of equivariant sheaves. In particular, it is shown in \cite{\bl} 
that the functor
$\F\mapsto q_*^G\F$ from the category of $G$-equivariant flat vector
bundles on $M$ to the category of flat vector bundles on $B=M/G$ is an
equivalence of categories. The inverse functor is the pull-back $q^*$
of vector bundles.

\subheading{2.6. Example} We will consider now the simplest situation 
leading to equivariant flat vector bundles. 

Let $\rho:G\to \End(V)$ be a linear representation of $G$ which is trivial
on the component $G_0$ of the unit element $e\in G$.  
Consider $\F=M\times V$ with the trivial connection and with
the diagonal $G$ action. It is an equivariant flat vector bundle. 
This bundle is trivial as a flat vector bundle but {\it it is not necessarily 
trivial as an equivariant flat vector bundle}. Thus, we may use the previous 
construction of pushforward to produce a flat bundle over $B=M/G$ starting 
from any linear representation of $G/G_0$.

Consider, for example, the case when $G=\ZZ_2$ acting in the standard way
$x\mapsto -x$ on the sphere $S^n$. Let $V=\CC$ be the unique non-trivial 
representation of $G$. Then the construction above produces the non-trivial 
flat line bundle over the real projective space $P^n$ (the M\"obius band).

\subheading{2.7. Equivariant cohomology twisted by an equivariant flat
bundle} Let $M$ be a $G$ manifold (the action is not supposed to be free),
and let $\F$ be an equivariant flat vector bundle over $M$. Our purpose now
is to define equivariant cohomology $H^\ast_G(M,\F)$. 

In a non-precise way we may describe our construction as follows. Let 
$EG\to BG$ be the universal principal bundle. Then we may consider the 
induced flat equivariant bundle $p^\ast\F$ over $EG\times M$, where 
$p:EG\times M\to M$ is the projection. Now, we want to form the pushforward
$q^G_\ast p^\ast\F$, where $q: EG\times M\to EG\times_G M =M_G$ 
is the projection; it is a flat vector bundle over the Borel's quotient
$M_G$. Now we define the equivariant cohomology $H^\ast_G(M,\F)$ as the 
cohomology of $M_G$ twisted by the flat vector bundle $q^G_\ast p^\ast\F$.

Unfortunately, we cannot really apply the construction described in the 
previous paragraph since our category is the category of smooth finite
dimensional manifolds 
and the universal principal bundle $EG\to BG$ is usually infinite 
dimensional. Instead, we will use finite dimensional approximations to the 
universal principal bundle. 

Consider any finite dimensional principal $G$-bundle $\pi: E\to B$ 
(with a right action of $G$ on $E$) over a smooth base
manifold $B$, and the corresponding mixing diagram
$$
    \CD
	E    @<<<     E\times M  @>p>>  M       \\
	@V{\pi}VV             @VVqV                 @VVV \\
	B  @<{\alpha}<<      E\times_G M  @>>>     M/G     .  
    \endCD
		\tag 23
$$
Here the action on the middle term 
$E\times M$ is the diagonal one $g(e,m)=(eg^{-1},gm)$ for all $g\in G$ and
$e\in E$, $m\in M$; it is a free action, since the action on $E$ is free; 
$q$ denotes the quotient map determined by this action of $G$.

Given a flat equivariant bundle $\F$ over $M$, we may first form the induced
bundle $p^\ast\F$ over $E\times M$ (cf. 2.3) and then apply the pushforward
construction (cf. 2.4) to obtain the flat bundle $q^G_\ast p^\ast\F$ over 
$E\times_G M$; the last flat bundle we will sometimes denote by $\F^E$ for 
short. Hence, the cohomology 
$$H^i(E\times_G M,\F^E)\tag24$$ 
is defined (as the cohomology of $E\times_G M$ with coefficients in the flat
bundle $\F^E$). We will show that for fixed $i$ this cohomology does not
depend on the choice of the principal bundle $\pi:E\to B$, provided that this 
principal bundle sufficiently closely approximates the universal bundle.

Let $n$ be an integer. A principal $G$-bundle $\pi:E\to B$ is called 
$n$-{\it acyclic} if the reduced complex cohomology $\tilde H^j(E,\CC)$
vanishes for all $j\le n$.

\proclaim{2.8. Lemma} Let $M$ be a $G$ manifold and let $\F$ be an equivariant
flat vector bundle over $M$. Let $E\to B$ and $E'\to B'$ be two $n$-acyclic 
principal $G$-bundles. Denote by $\F^{E}$ and $\F^{E'}$ the corresponding 
flat vector bundles on $E\times_GM$ and $E'\times_GM$ respectively, 
defined as explained above. Then
  $$
	H^i(E\times_GM,\F^{E})=H^i(E'\times_GM,\F^{E'})
		\quad \text{for all} \quad i=0,1\nek n-1.\tag25
  $$
\endproclaim
\demo{Proof}Set $E''=E\times E'$ and let $G$ act on $E''$ by $g:(e,e')\mapsto
(e\cdot g,e'\cdot g)$. Then $E''$ is an $n$-acyclic principal
$G$-bundle. It is enough to show that
$H^i(E''\times_GM,\F^{E''})=H^i(E\times_GM,\F^{E})$ for any
$i=0,1\nek n-1$. 

The natural projection $r:E''\times_GM\to E\times_GM$ defines a
locally trivial fibration whose fiber is isomorphic to $E'$. Also
$\F^{E''}= r^*\F^E$ as flat vector bundles. Hence, the restriction of 
$\F^{E''}$ on the fibers of $r$ is a trivial flat bundle.  Thus, the 
reduced cohomology $\tilde H^\ast(E''\times_G M,\F^{E''})$ may be calculated by 
means of the spectral sequence of the fibration $r$. The initial term 
of this spectral sequence is given by
$$
	E_2^{pq}\ =\ H^p(E\times_GM, \tilde\H^q(E',\CC)\otimes \F^E).\tag26
$$
Here $\tilde\H^q(E',\CC)$ is viewed as a local system of coefficients
over $E\times_GM$ determined by the fibration $r$, whose fiber over a 
point $x\in E\times_GM$ is equal to
the $q$-th reduced cohomology of the fiber $r^{-1}(x)\simeq E'$ and 
the product $\tilde\H^q(E',\CC)\otimes \F^E$ is understood as the 
tensor product of local coefficient systems.
 
The lemma follows now from the fact that $E'$ is $n$-acyclic.  
$\square$\enddemo

We can give now the main definition:

\proclaim{2.9. Definition} Let $M$ be a $G$ manifold and let $\F$ be an 
equivariant flat vector bundle over $M$. We define the twisted equivariant
cohomology
$$H^i_G(M,\F)\tag27$$
as the usual cohomology with coefficient in a flat vector bundle
$H^i(E\times_GM,\F^{E})$, where $E\to B$ is any $(i+1)$-acyclic 
principal $G$-bundle.

{\rm By Lemma 2.8 this definition does not depend on the choice of the
$(i+1)$-acyclic principal bundle $E$.}\endproclaim

\subheading{2.10} Note that the twisted equivariant cohomology
$H^\ast_G(M,\F)$ has a natural structure of a {\it graded module} over
the {\it graded equivariant cohomology ring} 
$H^\ast_G(M)\ =\ H^\ast_G(M,\CC)$. Also, the 
map $\alpha$ of the mixing
diagram (23) induces a homomorphism of graded rings
$$\alpha^\ast: H^\ast_G\ \to\ H^\ast_G(M)\tag 28$$
and so the twisted equivariant cohomology $H^\ast_G(M,\F)$ is naturally
a graded module over the cohomology ring 
$H^\ast_G\ =\ H^\ast_G(pt)=H^\ast(BG)$. It is well known that this
ring is a polynomial ring: 
$$
	H^\ast_G=\CC[u_1\nek u_l] \tag29
$$ 
with $l$ generators of even degree. Here $l$ is the rank of $G$,
i.e. the dimension of a maximal torus $T\subset G$. If $G=T$ is a
torus, these generators are all of dimension 2.

%----------------------------------------------------------------
\subheading{2.11. The Localization Theorem} One may formulate versions
of the Localization Theorem for twisted equivariant cohomology. We will 
mention one such statement here since we will need it in the present paper.
The proof repeats the well-known arguments (which can be found in \cite{Hs}
or \cite{AB2}) and so it will be skipped.

We will assume in this subsection that $G$ {\it is a torus}.

Let 
$$M^G=\{x\in M:\ g\cdot x=x \ \text{for any} \ g\in G\}$$ 
denote the set of fixed points of the action of $G$ on $M$. 
The inclusion $i:M^G\to M$ defines a homomorphism
$$
	i^*:H^\ast_G(M,\F)\to H^\ast_G(M^G,\F_{|_{M^G}})\tag30
$$
of $H^\ast_G$-modules.
The Localization Theorem states that {\it the kernel and cokernel of
$i^*$ are torsion modules. In particular,}
$$
	\rank H^\ast_G(M,\F)=\rank H^\ast_G(M^G,\F_{|_{M^G}})=
		\dim_\CC H^\ast(M^G,\F_{|_{M^G}}),\tag31
$$
where the rank is understood over the field of fractions of the cohomology 
ring $H^\ast_G$.
In the second equality we used that
$$
	H^\ast_G(M^G,\F_{|_{M^G}})\ =\ 
	H^\ast(M^G,\F_{|_{M^G}})\otimes_\CC H^\ast_G 
$$
as $H^\ast_G$-modules.

\subheading{2.12. Example} As an example 
consider the case of {\it the circle action} $G=S^1$. 
Applying the Localization Theorem  we obtain that the kernel and cokernel
of the map (30) are finitely generated torsion modules over $H_G^\ast = \CC[u]$ 
and so they 
are finite dimensional as vector spaces. It follows that 
$$
\dim_{\CC} H^i_G(M,\F) \ =\ \sum_{j \equiv i\mod 2} 
\dim H^j(M^G,\F_{|_{M^G}})\tag32
$$
for large values of the dimension $i$.

\heading \S 3. Some remarks about basic 1-forms  \endheading

In this section we will review some simple (and, probably, well-known)
properties of basic forms which will be used in the paper.

Let us recall that a 1-form $\tet\in \Ome^1(M,\RR)$ on a $G$ manifold
$M$ is called {\it basic} if it is $G$ invariant and its restriction
on any orbit of the action of $G$ vanishes.

\subheading{3.1} As the first remark let us mention the construction of {\it 
descent}. Suppose that the action of $G$ on $M$ is {\it free} 
and that $\tet$ is a
basic 1-form on $M$. Then (by the Gleason lemma, cf. [Hs], page 9) the
quotient space $B=M/G$ is a smooth manifold and the quotient map
$q:M\to B$ is a locally trivial fibration. As explained in  
\cite{KN, Ch. XII \S 1}, there exists a
unique 1-form $\otet$ on $B$ such that $q^*\otet=\tet$. The form
$\otet$ is closed if $\tet$ is closed. 

We will say that the form $\tet$ {\it descents} to $\otet$.

\subheading{3.2} Now we are going to mention the precise relation between the 
construction of descent and the pushforward construction of section 2.5.

Suppose again that $G$ acts freely on $M$ and let $q:M\to B$ denote the 
quotient map. Suppose that $\E$ is a line bundle over $B$ supplied with
two flat
connections $\nabla_1$ and $\nabla_2$. The difference $\nabla_1-\nabla_2$
is a closed 1-form on $B$. Consider the induced vector bundle $q^\ast\E$ 
over $M$; the flat connections $\nabla_1$ and $\nabla_2$ {\it determine 
uniquely the equivariant
flat connections} $\tilde \nabla_1$ and $\tilde \nabla_2$ on $q^\ast\E$,
respectively. Then, we claim that \newline
{\it the difference 
$\tilde \nabla_1\ -\ \tilde \nabla_2$ is a closed basic 1-form on $M$ which
descents to the form $\nabla_1 \ -\ \nabla_2$ on $B$.} 

Observe, that in the above situation the pushforward of the flat equivariant
bundle $(q^\ast\E,\tilde\nabla_\nu)$ is the flat bundle $(\E,\nabla_\nu)$,
where $\nu\ =\ 1, 2$.

Next we will consider the question which cohomology
classes in $H^1(M,\RR)$ may be represented by a closed basic 1-form. Here
we will not assume that $G$ acts freely.

\proclaim{3.3. Lemma} Let $M$ be a $G$-manifold, where $G$ is a compact
Lie group. A cohomology class $\xi\in H^1(M,\RR)$ can be
represented by a closed basic 1-form if and only if it belongs to
the image of the homomorphism 
$$H^1_G(M,\RR)\to H^1(M,\RR).\tag33$$
\endproclaim
\demo{Proof} Let $\pi:E\to B$ be a smooth $2$-acyclic (cf. 2.7) principal 
$G$-bundle. The space
$E\times_GM$ has the structure of a locally trivial fibration
over the  space $B$ with fiber $M$. 
The inclusion of the fiber $i:M\to E\times_GM$ induces
the homomorphism
$$
	i^*:H^1(E\times_GM)=H^1_G(M)\to H^1(M)\tag34
$$
which can be identified with the homomorphism (33).
We will denote by  $p:E\times M\to M$ the projection.  

If $\theta$ a closed basic 1-form on $M$ then $p^\ast\theta$ is a closed  
basic 1-form on $E\times M$; by the construction of descent (cf. 3.1)
it determines a closed 1-form $\otet$ on $E\times_G M$ such that
$q^\ast\otet\ =\ p^\ast\theta$. Clearly, $i^\ast\otet\ =\ \theta$.
Thus, we have $[\otet]\in H^1(E\times_G M)\ =\ H^1_G(M)$ and 
$i^\ast[\otet]\ =\ [\theta]$, so $[\theta]$ belongs to the image of (34).

Conversely, if a cohomology class $\xi\in H^1(M,\RR)$ belongs to the image
of (33) then there exists $\oxi\in H^1(E\times_G M)$ with $i^\ast\oxi=\xi$.
We may realize $\oxi$ by a closed 1-form $\otet$; 
then $\theta=i^\ast \otet$ is a basic 1-form on $M$ realizing $\xi$.
$\square$. \enddemo
  
We finish this section with the following useful lemma, which shows that 
in many applications any $G$-invariant 1-form is automatically basic.

\proclaim{3.4. Lemma} Let $M$ be a {\it connected} $G$-manifold. Assume
that the  fixed points set of the action of $G$ on $M$ is not
empty. Then any closed $G$-invariant 1-form on $M$ is basic.
\endproclaim
\demo{Proof}Let $\w$ be a closed $G$-invariant 1-form on $M$. We have 
to show that  
$$
	\iot(X_M)\w=0           \tag35
$$ 
for any $X\in \grg$. Clearly, it is enough to prove (35) for those
$X\in \grg$ for which the one parameter subgroup 
$\{ \exp(tX)| \ t\in \RR\}$ generated by $X\in \grg$, is compact. 

Fix $X\in \grg$ such that the subgroup $\{ \exp(tX)| \ t\in \RR\}$
of $G$ is compact. For each $x\in M$ we denote by
  $$
	\Gam_x = \{ \exp(tX)\cdot x| \ t\in \RR\}\tag36
  $$
  its orbit under the action of $\{ \exp(tX)| \ t\in \RR\}$. Since $\w$
  is $G$-invariant, (35) is equivalent to
  $$
	\int_{\Gam_x}\w = 0 \quad \text{for any} \quad X\in \grg, \quad x\in M.
			\tag37
  $$
The fixed points set of the action $\{ \exp(tX)| \ t\in \RR\}$ on
$M$ is supposed to be not empty. Hence, each orbit $\Gam_x$ is free
homotopic to a point. Since $\w$ is closed, this implies (37).
$\square$\enddemo

\heading \S 4. Proof of Theorem 1.7\endheading

Throughout this section we will assume that we are
in the situation of Theorem 1.7, i.e. we have a closed $G$ manifold $M$, a
flat equivariant vector bundle $\F$ over $M$ and a closed basic 1-form
$\theta$ on $M$ which is non-degenerate in the sense of Bott. We will denote
by $\xi\in H^1(M,\RR)$ the cohomology class of $\theta$. As shown in 2.2,
the form $\theta$ determines a flat equivariant line bundle $\E_\theta$
over $M$; in fact, we will consider a one-parameter family of flat
equivariant bundles $\F\otimes \E_{t\theta}$ over $M$, where $t\in\RR$ is 
a real parameter.

First we are going to obtain an approximate version of Theorem 1.7 (cf. 
Proposition 4.4). 

Consider a principal smooth $G$ bundle $E\to B$ over a smooth manifold B;
here $G$ acts on $E$ from the right. As in subsections 2.7 - 2.9, 
we will consider the family of flat bundles 
$$(\F\otimes \E_{t\theta})^E\ 
=\ \F^E\otimes \E_{t\theta}^E,\quad t\in\RR\tag38$$ 
over $E\times_G M$, which are obtained form the family of flat equivariant 
bundles 
$\F\otimes \E_{t\theta}$ over $M$ by first inducing flat bundles over
$E\times M$ and then performing the pushforward construction 2.5 with respect
to the projection $E\times M\to E\times_G M$.

\proclaim{4.1. Lemma-Definition} Consider the dimension of the cohomology 
$$d(t)\ =\ \dim H^i(E\times_G M, \F^E\otimes \E_{t\theta}^E)\tag39$$ 
as a function of $t\in \RR$.
Then there exists a finite subset $S\subset \RR$ such that $d(t)=d_0=\const$
for all $t\notin S$ and $d(t)>d_0$ for $t\in S$. We will call $d_0$ the 
$i$-dimensional $E$-equivariant Novikov number and will denote it by 
$\bet_i(\xi,\F;E)$; this number depends only on the cohomology class $\xi$
of $\theta$ by Lemma 2.3. \endproclaim

\demo{Proof} By 3.2 there exists a closed 1-form $\tilde\theta$ on
$E\times_G M$ such that for all $t$ the bundle $\E^E_{t\theta}$ coincides
(as a flat bundle) with the flat line bundle over $E\times_G M$ with the 
connection $d + t\tilde\theta\wedge\cdot$. The Lemma follows now from 
Lemma 1.3 of \cite{BF1}. $\square$ 
\enddemo

\subheading{4.2} Observe that by Lemma 3.3, any class $\xi$ in the 
image of the map
$H^1_G(M,\RR)\to H^1(M,\RR)$ can be realized by a closed basic form and  
so the $E$-equivariant Novikov numbers $\bet_i(\xi,\F;E)$ are defined for 
all such classes.

The {\it $E$-equivariant Novikov polynomial} $\calN_{\xi,\F}^E(\lam)$
associated to a cohomology class $\xi\in \im[H^1_G(M,\RR)\to H^1(M,\RR)]$ 
and to a $G$-equivariant flat vector bundle $\F$ is defined by the
formula 
$$
	\calN_{\xi,\F}^E(\lam)=\sum_{i=0}^k\lam^i\bet_i(\xi,\F;E),
			\qquad k=\dim_{\CC} (E\times_GM).\tag40
$$

\subheading{4.3. The $E$-equivariant Morse polynomial} 
Let $Z$ be a connected component of the critical set $C$ of $\w$ and
let $G_Z=\{g\in G: \ g\cdot Z\subset Z \}$
denote the stabilizer of the component $Z$. Let $|G:G_Z|$ denote the
index of $G_Z$ as a subgroup of $G$. Note that $G_Z$ contains
the connected component of the unity in $G$ (in particular, if
$G$ is connected then $G_Z=G$). Hence $|G:G_Z|$ is finite.

Recall that in Section 1.5 we associated to each connected component $Z$
the flat line bundle $o(Z)$ over $Z$ and an integer $\ind(Z)$. Let $\F_{|_Z}$ denote the restriction
of $\F$ to $Z$. Clearly, $o(Z)$ and $\F_{|_Z}$ are $G_Z$-equivariant
flat bundles.

The principal bundle $E$ may be considered as a principal $G_Z$
bundle over $E/G_Z$. Applying the induction (cf. 2.4) and then 
pushforward (cf. 2.5) to the flat equivariant bundles $\F_{|_Z}$ and $o(Z)$ 
we obtain flat vector bundles $\F_{|_Z}^E$
and $o(Z)^E$ over $Z^E=E\times_{G_Z} M$. The numbers
$$
	\bet^E_i(Z,\fo)=
		\dim_{\CC}H^i(Z^E,\F^E_{|_{Z}}\otimes o(Z)^E),
			\qquad i=0,1\nek \dim Z^E\tag41
$$
are called the {\it  $E$-equivariant Betti numbers} of $Z$.  
Here $H^i(Z^E,\F^E_{|_{Z}}\otimes o(Z)^E)$ denotes the cohomology 
of $Z^E$ with coefficients in the flat vector bundle
$\F^E_{|_{Z^E}}\otimes o(Z)^E$.

Consider now the {\it twisted $E$-equivariant Poincar\'e
polynomial} of $Z$
$$
	\calP_{Z,\F}^E(\lam)=\sum \lam^i \bet^E_i(Z,\fo)\tag42
$$ 
and define using it the following {\it $E$-equivariant Morse counting
polynomial} of $\w$
$$
	\calM_{\w,\F}^E(\lam)=
		\sum_Z \lam^{\ind(Z)}|G:G_Z|^{-1}\calP_{Z,\F}^E(\lam),\tag43
$$
where the sum is taken over all connected components $Z$ of $C$.

The following result can be considered as an approximate version of 
Theorem 1.7.
%-----------------------------
\proclaim{4.4. Proposition}There exists a polynomial $\calQ(\lam)$ with
non-negative integer coefficients, such that
   $$
	\calM_{\w,\F}^E(\lam)- \calN_{\xi,\F}^E(\lam)=
		(1+\lam)\calQ(\lam).\tag44
   $$
\endproclaim
%-----------------------
\demo{Proof} Let $p:E\times M\to M$ denote the projection and let 
$q:E\times M\to E\times_G M$ denote the quotient map (as in the mixing 
diagram (23)). It is shown in section 3.1 that there exists a smooth closed
1-form $\theta^E$ on $E\times_G M$ such that $p^\ast\theta = q^\ast\theta^E$.
{}From the construction of $\w^E$ it is clear that it is non-degenerate in 
the sense of Bott assuming that $\theta$ is non-degenerate.

The proof of Proposition 4.4 is based on the non-equivariant Novikov-Bott 
inequalities, established by the authors in \cite{\bfa, Theorem 0.3} and in
\cite{\bfb, Theorem 4}, applied to the form $\w^E$. We shall explain the
details.

Let $C$ denote the critical set of $\w$. Then the critical set
$C^E$ of $\w^E$ coincides with $E\times_G C\subset E\times_G M$. If $Z$
   is a connected component of $C$, then 
   $$
		Z^E=E\times_{G_Z} Z\subset E\times_G M\tag45
   $$ 
   is a connected component of $C^E$. 

   Recall that in Section 1.5 we associated to a connected
   component $Z$ of $C$ the bundle $o(Z)$ and the number
   $\ind(Z)$. Let $o(Z^E)$ and $\ind(Z^E)$ be the corresponding
   objects associated to the component $Z^E$ of $C^E$.  Simple
   calculations (cf. \cite{AB1, Proposition 1.5}) 
show that $\ind(Z)=\ind(Z^E)$ and $o(Z^E)=o(Z)^E$. It
   follows that the {\it $E$-equivariant Betti numbers}
   $\bet^E_i(Z,\fo)$ of $Z$ may be calculated by the formula
   $$
	\bet^E_i(Z,\fo)=
		\dim_{\CC}H^i(Z^E,{\F^E}_{|_{Z^E}}\otimes o(Z^E)),
			\qquad i=0,1\nek \dim Z^E.\tag46
   $$
   Here ${\F^E}_{|_{Z^E}}$ is the restriction of $\F^E$ on $Z^E$ and
   $H^i(Z^E,{\F^E}_{|_{Z^E}}\otimes o(Z^E))$ denote the cohomology of
   $Z^E$ with coefficients in the flat vector bundle
   ${\F^E}_{|_{Z^E}}\otimes o(Z^E)\ 
   =\ q_\ast^G p^\ast(\F_{|_Z}\otimes o(Z))$.

   The correspondence $Z\mapsto Z^E=Z\times_{G_Z}E$ is a surjection
   from the set of connected components of $C$ to the set of connected
   components of $C^E$. Moreover the preimage of $Z^E$ contains
   exactly $|G:G_Z|$ connected components of $C$. Hence,
   $$
	\calM^E_{\w,\F}(\lam)=
	   \sum_{Z}\lam^{\ind(Z)} \cdot  \sum_{i=0}^{\dim Z^E}
	    \lam^i \dim_{\CC} H^i(Z^E,{\F^E}_{|_{Z^E}}\otimes o(Z^E)).  \tag47
   $$ 
   This means that the $E$-equivariant Morse counting polynomial of
   $\w$ coincides with the Morse counting polynomial of $\w^E$ as it
   is defined in \cite{\bfa,\bfb}:
   $$\calM^E_{\theta,\F}(\lam) = \calM_{\theta^E,\F^E}(\lam).$$
   The $E$-equivariant Novikov
   polynomial $\calN^E_{\xi,\F}(\lam)$  by definition is equal to the
   Novikov polynomial associated to the cohomology class of $\w^E$ in
   $H^1(M^E,\RR)$ (cf. \cite{\bfa,\bfb}). The proposition follows now
   from Theorem 0.3 of \cite{\bfa} (see also \cite{\bfb, Theorem 4}).
$\square$\enddemo

\subheading{4.5. Proof of Theorem 1.7} Let
$$\calM_{\w,\F}^G(\lam)\ =\ \sum \alpha_i\lambda^i,\qquad 
\calN_{\xi,\F}^G(\lam)\ =\ \sum \beta_i\lambda^i\tag48$$
and denote $\gamma_i=\alpha_i-\beta_i$. Then Theorem 1.7 states that
$$\gamma_p-\gamma_{p-1}+\gamma_{p-2}-\dots +(-1)^p\gamma_0\ge 0\tag49$$
for all $p=0,1,2,\dots$. For any $p$ this inequality involves only 
equivariant cohomology (of $M$ and of the critical manifolds) of dimension
less than $p+1$. Hence, each of these inequalities
follows from Proposition 4.4 and from the definition of twisted equivariant
cohomology (cf. Definition 2.9). This completes the proof. $\square$

\heading{\S 5. Application: fixed points of a symplectic torus action}
\endheading

In this section we apply Theorem 1.7 to study the topology of the fixed point set of a
symplectic torus action. We will see that inequalities given by Theorem 1.7 are exact
in this case. This gives new  relations between the
Betti numbers of the fixed points set of a torus action on a symplectic
manifold, and the equivariant Novikov numbers.

%----------------------
\subheading{5.1} Assume that the $n$-dimensional torus
$T=S^1\times\cdots\times S^1$ acts by symplectomorphisms on a compact
connected symplectic manifold $(M,\ome)$. Denote by
$$
	M^T\ =\ \big\{x\in M| \ t\cdot x=x \ \ 
		\text{for any} \ \ t\in T\big\} \tag50
$$
the fixed point set of the action.

Let $\grt$ be the Lie algebra of $T$. Choose $X\in \grt$ such
that the one-parameter subgroup \/ $\{ \exp(sX)| \ s\in\RR\}$ \/
generated by $X$ is dense in $T$. Let $X_M$ be the vector field on
$M$ defined by the infinitesimal action of $\grt$ on $M$. Let
$\iot(X_M)$ denote the interior multiplication by $X_M$.
Set
$$
	\w=\iot(X_M)\ome. \tag51
$$
Then $\w$ is a closed $T$-invariant 1-form on $M$.
We will call $\theta$ {\it the generalized moment map} of the torus action.

McDuff \cite{\mcd} constructed examples showing that this form may be not 
cohomologous to zero.
Clearly, the critical set of
$\w$ coincides with the set $M^T$ of fixed points.
Hence, one can use Theorem 1.7, to estimate the Betti numbers
of $M^T$. It is shown in \cite{\au}, cf. \S\S 2.1, 2.2, that $\w$ is
non-degenerate in the sense of Bott, each connected component $Z$ of
$M^T$ has even index, and the flat line bundle $o(Z)$ is trivial.

In order to satisfy the assumption of Theorem 1.7 that $\theta$ is basic, we will assume
that either the fixed point set $M^T$ is not empty (then $\theta$ is basic by Lemma 3.4),
or that $T$ is the circle (and then $\theta$ is clearly basic).

%---------------------------------------
\subheading{5.2}
Let $\F$ be a $T$-equivariant flat vector bundle on $M$ (for instance,
one may always take the trivial line bundle for $\F$.) 
For any connected component $Z$ of $M^T$ we denote by 
$$
	\calP_{Z,\F}(\lam)=
		\sum_{i=0}^{\dim Z}\lam^i \dim_{\CC} H^i(Z,\F_{|_Z})\tag52
$$
the Poincar\'e polynomial of $Z$ and by $\ind(Z)$ the index of $Z$ 
considered as a critical submanifold of the form $\w$. $\calP^T(\lam)$ 
will denote the $T$-equivariant Poincar\'e series of a 
point (as in 1.9). It is well known that 
$\calP^T(\lam)= (1-\lam^2)^{-n}.$
The equivariant Morse polynomial of $\w$ is given by the formula
$$
	\calM_{\w,\F}^T(\lambda)= 
	   (1-\lam^2)^{-n} \cdot
		\sum_Z \lam^{\ind(Z)}\calP_{Z,\F}(\lam)
		\tag 53
$$
where the sum is taken over all connected components of the fixed
points set $M^T$.

Let $\xi=[\w]\in H^1(M,\RR)$ denote the cohomology class of $\w$ and
let  $\calN_{\w,\F}(\lam)$ denote the equivariant Novikov series of $\F$
associated to this cohomology class, cf. 1.4.
\proclaim{5.3. Theorem} In the situation described in 5.1 and 5.2, assume that either
the fixed point set $M^T$ is not empty, or $T$ is the circle.
Then the following identity holds
$$(1-\lam^2)^{-n}\cdot \sum_Z \lam^{\ind(Z)}\calP_{Z,\calF}(\lam) =
\calN_{\xi,\calF}^T(\lam),\tag54$$
where the sum on the left is taken over all connected components $Z$
of the fixed point set $M^T$.
\endproclaim

This theorem gives new relations between the homology of the
fixed point set $M^T$ and the equivariant Novikov numbers. Compare
\cite{G}, \cite{Hs}.

\demo{Proof} The following arguments generalize the result of V. Ginzburg
\cite{Gi}; cf. also \cite{AB2}, page 26. Note that papers \cite{Gi}
and \cite{AB2} deal with Hamiltonian actions, i.e. assuming that the cohomology class of
$\theta$ is zero.

Choose a generic $t$. More precisely, using Lemma 1.3, we may find a 
countable subset $S\subset \RR$ such that for $t\notin S$ the dimension
of $H^i_G(M,\F\otimes \E_{t\theta})$ equals the Novikov number
$\beta_i^T(\xi,\F)$ for all $i$. We will assume in the sequel that $t\notin S$.

Recall from section 2.10, that the twisted equivariant cohomology
$$H^\ast_T(M,\calF\otimes \calE_{t\theta}) =
\bigoplus_{i=0}^{\infty} H^i_T(M,\calF\otimes \calE_{t\theta})$$
has a natural structure of a graded module over the ring $H^\ast_T = \CC[u_1,\dots,u_n]$.
Similarly, each term of the spectral sequence
$E_r^{pq}$ (where $r = 2,3,\dots$) of the fibration
$$
\CD
M @>i>> ET\times_T M\\
&  &       @V{\pi}VV\\
   &   &      BT,
\endCD
$$
with coefficients in the equivariant flat vector bundle $\calF \otimes\calE_{t\theta}$,
is a graded module over $H^\ast_T$ (cf. notation (29)).
Here the spectral sequence $E_r^{\ast,\ast}$
is considered with respect to
the total grading: $E_r = \oplus E_r^n$, where
$$E^n_r = \bigoplus_{p+q=n} E^{pq}_r.$$
Our aim is to show that all the differentials of this spectral sequence vanish. Note that
from the
general properties of the spectral sequence of a fibration we know that the differentials are
$H^\ast_T$-module homomorphisms.

We will deal with the {\it ranks} of the terms of this
spectral sequence, which are defined as the dimensions over the field of rational functions
$$\tilde H^\ast_T = \CC(u_1,\dots,u_n)$$
of their localizations $\tilde H^\ast_T\otimes_{H^\ast_T} E_r$; we will denote these
ranks by $\rank E_r$.
Note that the initial term
$$E_2 = \bigoplus_{p,q} \ H^p_T\otimes H^q(M,\calF\otimes \calE_{t\theta})$$
of this spectral sequence {\it is free}
as $H^\ast_T$-module and its rank is
$$\rank E_2 = \sum_{i=0}^{\dim M} \dim H^i(M,\calF\otimes\calE_{t\theta}).$$
On the other hand, the rank of the limit term $E_{\infty}$ equals to
$$\rank E_\infty = \sum_{i=0}^{\dim M^T} \dim H^i(M^T, \calF|_{M^T}),\tag55$$
as follows from the Localization theorem, cf. 2.11. Since $\rank E_2 \ge \rank E_\infty$,
we obtain the inequality
$$\sum_{i=0}^{\dim M} \dim H^i(M,\calF\otimes\calE_{t\theta}) \ge
\sum_{i=0}^{\dim M^T} \dim H^i(M^T, \calF|_{M^T}).\tag56$$
On the other hand, the non-equivariant Novikov - Bott inequalities of \cite{BF1}, Theorem
0.3 (cf. also \cite{BF2}, Theorem 4) give (by substituting $\lam =1$) the opposite inequality
to (56).

Thus we obtain that the rank of the terms of the spectral sequence is constant
$\rank E_r = const$. Since the initial term is free, the first nontrivial differential would
reduce the rank; therefore all the differentials vanish.

As another conclusion of the above arguments we obtain the identity
$$\sum_Z \lam^{\ind Z}(\sum_{i=0}^{\dim Z}\lam^i \dim H^i(Z,\calF|_Z)) =
\sum_{i=0}^{\dim M} \lam^i\dim H^i(M,\calF\otimes\calE_{t\theta}),\tag57$$
which follows from Theorem 0.3 of \cite{BF1}. In fact, we know that the polynomial $\calQ(\lam)$
in Theorem 0.3 of \cite{BF1} has nonnegative coefficients, and also we know that
$\calQ(1)=0$; therefore
we obtain that $\calQ(\lam)$ is identically zero.

Returning again to the spectral sequence, we obtain (by comparing the Poincar\'e
power series of the initial term with the Poincar\'e series of the limit) that
$$(1-\lam^2)^{-n} \cdot\sum_{i=0}^{\dim M} \lam^i\dim H^i(M,\calF\otimes\calE_{t\theta})
= \sum_{i=0}^\infty \lam^i\dim H^i_T(M, \calF\otimes \calE_{t\theta}).\tag58$$
This, combined with (57), gives the required identity (54). \qed

\enddemo

\subheading{5.4. Remark} 
Note that a similar statement concerning {\it perfectness or exactness} of our
equivariant Novikov inequalities (10)
can be proven in a more general situation of Theorem 1.7 assuming only that
the set of critical points of a basic 1-form coincides with the fixed point set of an
action of a compact connected Lie group. All the above arguments can then be applied.

Next we present some easy corollaries of Theorem 5.3.

\proclaim{5.5. Corollary} Let $M$ be a symplectic manifold with a symplectic circle action.
$M$ has no fixed points
if and only if all the equivariant Novikov numbers $\beta^T_i(\xi,\calF)$,
where $i=0,1,2\dots$,
vanish.
Here $\calF$ denotes the trivial flat line bundle over $M$, and $\xi=[\theta]$ is the
cohomology class of the generalized moment map $\theta$ (cf. (51)).
\endproclaim

\demo{Proof} First we should mention that the form $\theta$ (defined by (51)) is clearly
basic and so the equivariant Novikov numbers are well defined. Now one simply applies
Theorem 5.3. \qed
\enddemo

Here is another simple corollary:

\proclaim{5.6. Theorem} Let $M$ be a symplectic manifold with a
symplectic action of the oriented circle $T=S^1$, and let $\F$ be a flat equivariant
bundle of rank $d$ over $M$. Suppose that the fixed point set $M^T$
is non-empty and all the odd-dimensional cohomology 
$H^{odd}(M^T,\F_{|_{M^T}})$ 
vanishes (note that this is automatically true if $M^T$
consists of isolated points!) Fix a positive vector $X\in\grt$ and 
consider the 1-form $\theta=i(X_M)\omega$ as in (51).
Then: 
\roster

\item all odd-dimensional equivariant Novikov numbers $\beta_{2i-1}^T(\xi,\F)$ 
vanish;
\item the even dimensional equivariant Novikov numbers increase
$$\beta_{0}^T(\xi,\F) \le \beta_{2}^T(\xi,\F) \le \beta_{4}^T(\xi,\F)\le 
\dots,\tag59$$
and they stabilize 
$$\beta_{2i}^T(\xi,\F) = \beta_{n}^T(\xi,\F) = d\cdot\chi(M)\tag60$$
for $i\ge n/2$, where $n = \dim M$;
\item assume additionally that the set of fixed points $M^T$ is finite.
Let $m_i$ denote the number of fixed points of $T$, which have index $i$
as critical points of $\theta$. Then
$$m_{2i} = d^{-1}[\beta_{2i}^T(\xi,\F) - \beta_{2i-2}^T(\xi,\F)]
\quad\text{and}\quad m_{2i-1} = 0 \tag61$$
for all $i$. Moreover, the following symmetry relation
$$m_i = m_{n-i}\tag62$$
holds for all $i$. The total number of fixed points of $T$ equals to
$$d^{-1}\beta_n^T(\xi,\F) = \chi(M).\tag63$$

\endroster
\endproclaim

\demo{Proof} Applying Theorem 5.3 we obtain
$$(1-\lam^2)^{-1}\sum_Z \lam^{\ind(Z)}\calP_{Z,\F}(\lam) 
= \calN_{\xi,\F}^T(\lam).\tag64$$
Assuming that $H^{odd}(M^T,\F|_{M^T})=0$,  we see that the power series 
on the left side of (64) has only even powers of $\lam$.
This implies that the Novikov
power series $\calN_{\xi,\F}^T(\lam) $ has also only even powers of $\lam$.

Formula (59) and the fact that the Novikov numbers stabilize (i.e.
$\beta_{2i}^T(\xi,\F) =$ $\beta_{n}^T(\xi,\F)$ for $i\ge n/2$) follow from
(64) since the left hand side of (64) is clearly a rational function of the form
$(1-\lam^2)^{-1}p(\lam)$, where the polynomial $p(\lam)$ is of degree $\le n$ and has
only even powers with nonnegative integral coefficients. Expanding it into a power series,
we see that the {\it stable Novikov number} $\beta_{n}^T(\xi,\F)$ equals 
$p(1)$, which is the same as $d\cdot\chi(M)$ (by corollary 0.4 of \cite{BF1} and since
all the indices $\ind (Z)$ are even).

If the fixed point set $M^T$ consists of isolated points, then the above 
polynomial $p(\lam)$ is just
$\sum m_i\lam^i$ and the formula (61) follows from (64) 
by comparing the coefficients. 
Summing up formulas (61) gives (63).  Formula (62) follows by reversing 
the orientation of the circle and using  
the obvious relation $\beta_i^T(\xi,\F) = \beta_i^T(- \xi,\F)$ between the 
equivariant Novikov numbers.
\qed
\enddemo

Similar statement holds for the torus actions as well.

Note also that in  the K\"ahler case, $\xi=0$ by \cite{Fr} and, thus, the equivariant
Novikov numbers in Theorem 5.6 can be replaced by the dimensions of the
equivariant 
cohomology in the corresponding dimension. This observation shows that in this
special case our Theorem 5.6 gives the following well-known statement 
(cf. \cite{AB2}, page 23):

\proclaim{5.7. Corollary} Suppose that an oriented circle $T$ acts 
holomorphically 
on a K\"ahler manifold $M$ of dimension $n$
such that the K\"ahler form is preserved and such that the fixed point set
is non-empty and consists of finitely many isolated points.
Then for odd $i$ the equivariant cohomology $H^i_T(M)$ vanishes;
for even numbers $i$ the dimensions of the equivariant cohomology 
$H^i_T(M)$ increase
$$\dim H^0_T(M)\le \dim H^2_T(M)\le \dim H^4_T(M)\le \dots\le
\dim H^n_T(M)\tag65$$
and stabilize 
$$\dim H^i_T(M)=\dim H^n_T(M)\tag66$$
for even $i\ge n = \dim M$.  Additionally, for $i$ even, $0\le i\le n$, 
the difference
$$\dim H^i_T(M) - \dim H^{i-2}_T(M)$$
equals to 
$$m_i = \dim H^i(M),\tag67$$ 
where $m_i$ denotes the number of fixed points of 
the action of the circle $T$ having index $i$.
The total number of fixed points of $T$ is given by the dimension of the 
equivariant cohomology
$$\dim_{\CC}H^n_T(M) = \chi(M).\tag68$$
\endproclaim

\demo{Proof} As we noticed above, in the situation of Corollary 5.7 
$\xi=0$ by \cite{Fr},
and so there is a Hamiltonian $f:M\to\RR$ such that the form $\theta$
(given by (51)) is $\theta=df$. Then $f$ is a Morse function with critical
points of even indices only. The Morse theory implies that $f$ is perfect 
and this gives formula (67). The other statements of Corollary follow
directly from Theorem 5.5. \qed
\enddemo

Observe, that applying the Morse theory to the Hamiltonian $f$ (as above)
we find that the manifold $M$ in Corollary 5.7 is necessarily simply
connected. That is why we cannot expect to gain additional information
by considering general equivariant flat vector bundles $\F$.

Corollary 5.7 is valid in a more general situation of Hamiltonian circle
actions with isolated fixed points -- we have only used the 
condition $\xi=0$ in the proof (cf. also \cite{AB2}, page 23). 
Recall that by result of D. McDuff \cite{McD} 
any symplectic circle action on a 4-dimensional manifold is Hamiltonian if
it has fixed points.
Also, the Corollary remains true for  Lefschetz manifolds
(which include the K\"ahler manifolds),
i.e. for symplectic manifolds $M$ such that the multiplication by 
$\omega^{n-1}$
defines an isomorphism $H^1(M,\RR)\to H^{2n-1}(M,\RR)$. It was shown
by K.Ono \cite{O} that a symplectic action of a circle on a Lefschetz
manifold is Hamiltonian iff it has fixed points, cf. also \cite{McD}.

%To illustrate Corollary 5.6, consider the following example.
%Let $M$ be the standard unit sphere $S^2\subset \RR^3$ with center at 
%the origin, and let the circle $T$ act by rotations around the $z$ axis.
%Then the dimensions of the equivariant
%cohomology $H^\ast_T(M)$ are
%$$1, 0, 2, 0, 2, 0, 2, \dots$$
%and we have two fixed points (the North and the South poles) which have
%indices 0 and 2. Thus, $m_0=1$ and $m_2=1$ in agreement with Corollary 5.6.

E.Witten in \cite{Wi} constructed Morse type inequalities
which estimate the topology of the fixed point set of a circle action on a K\"ahler
manifold. He also assumes that the action is holomorphic and preserving the 
K\"ahler form. Cf. also \cite{MW}. A generalization of these inequalities of 
Witten for compact Lie groups was
announced recently in \cite{Wu}. S.Wu assumes however that the fixed points are 
isolated. 

It would be interesting to compare 
the results of this paper with the information given by \cite{Wi}, \cite{Wu}.

%------------------------------------------------------------------
%------------------------------------------------------------------

\enddocument